\documentclass[10pt]{amsart}

\usepackage{a4wide,amsmath,amssymb}

\newtheorem{thm}{Theorem}[section]
\newtheorem{lemma}[thm]{Lemma}
\newtheorem{prop}[thm]{Proposition}
\newtheorem{cor}[thm]{Corollary}
\theoremstyle{remark}
\newtheorem{rem}[thm]{Remark}

\theoremstyle{definition}

\newtheoremstyle{Claim}{}{}{\itshape}{}{\itshape\bfseries}{:}{ }{#1}
\theoremstyle{Claim}

\newenvironment{sistema}
{\left\lbrace\begin{array}{@{}l@{}}}
{\end{array}\right.}

\newcommand{\Z}{{\mathbb{Z}}}

\newcommand{\R}{\mathbb{R}}

\newcommand{\T}{\mathbb{T}}

\newcommand{\ep}{\epsilon}

\newcommand{\al}{\alpha}

\newcommand{\de}{\delta}

\newcommand{\la}{\lambda}

\theoremstyle{plain}

\def\sideremark#1{\ifvmode\leavevmode\fi\vadjust{
\vbox to0pt{\hbox to 0pt{\hskip\hsize\hskip1em
\vbox{\hsize3cm\tiny\raggedright\pretolerance10000
\noindent #1\hfill}\hss}\vbox to8pt{\vfil}\vss}}}

\begin{document}

\title[Asymptotic properties of $H$ with drift]{Asymptotic properties of\\ 
non-coercive Hamiltonians with drift}

\author{Martino Bardi}

\date{January 16, 2024}
 \thanks{
 The author is member of the Gruppo Nazionale per l'Analisi Matematica, la Probabilit\`a e le loro Applicazioni (GNAMPA) of the Istituto Nazionale di Alta Matematica (INdAM).
 He is partially supported by the King Abdullah University of Science and Technology (KAUST) project CRG2021-4674  ``Mean-Field Games: models, theory, and computational aspects", and by the project funded by the EuropeanUnion-NextGenerationEU under the National Recovery and Resilience Plan (NRRP), Mission 4 Component 2 Investment 1.1 - Call PRIN 2022 No. 104 of February 2, 2022 of Italian Ministry of University and Research; Project 2022W58BJ5 (subject area: PE - Physical Sciences and Engineering)  ``PDEs and optimal control methods in mean field games, population dynamics and multi-agent models"}
 
 \address{Department of Mathematics ``T. Levi-Civita", University of Padova, Via Trieste 63, 35121 Padova, Italy} \email{bardi@math.unipd.it}

\begin{abstract} 
We consider Hamiltonians associated to optimal control problems for affine systems on the torus. They are not coercive and are possibly unbounded from below in the direction of the drift of the system. The main assumption is the strong bracket generation condition on the vector fields. We first prove the existence of a critical value of the Hamiltonian by means of the ergodic approximation. Next we prove the existence of a possibly discontinuous viscosity solution to the critical equation. We show that the long-time behaviour  of solutions to the evolutive Hamilton-Jacobi equation is described in terms of the critical constant and a critical solution. As in the classical weak KAM theory we find a fixed point of the Hamilton-Jacobi-Lax-Oleinik semigroup, although possibly discontinuous. Finally we apply the existence and properties of the critical value to the periodic homogenization of stationary and evolutive H-J equations.
\end{abstract}
\maketitle


\section{Introduction}
\label{intro}
This paper gives 
 several asymptotic results about non-coercive Hamiltonians $H : \T^n\times \R^n \to \R$ possibly unbounded from below. By coercive we mean
 \[
\lim_{|p|\to\infty} \inf_{x\in \T^n} H(x,p)=+\infty .
\]
We are interested in Hamiltonians of the form
\begin{equation}
\label{Ham1}
H(x,p)= b(x)\cdot p + \sup_{a\in A} \left\{ p \cdot F(x) a - L(x, a) 
\right\},
\end{equation}
where $A\subseteq \R^m$, possibly all $\R^m$, $b$ is Lipschitz, $L$ is continuous, bounded from below, and quadratic in $a$ at infinity, $F(x)$ is a $n\times m$ matrix, at least Lipschitz, and whose columns $f^i$ satisfy a Lie-bracket rank condition. All data are $\Z^n$-periodic in $x$. Note that $H$ is unbounded from below if $b(x)\notin \text{span} \{f^1(x),..., f^m(x)\}$ for some $x\in \T^n$.

Let us list the properties we are going to study and briefly review what is already known, mostly in the coercive case and often under additional smoothness and convexity conditions on $H$.

{\bf 1.} The {\em critical} or {\em ergodic} equation associated to $H$ is
\begin{equation}
\label{eq:crit1}
\lambda+ H(x, Du
 )=0 
 \quad\text{ in } \R^n ,
\end{equation}
in the unknowns $\lambda\in \R$ and $u\in C(\T^n)$. The {\em critical value} of $H$ is a number $\bar\lambda\in \R$ such that \eqref{eq:crit1}  has no viscosity subsolutions for $\lambda > \bar\lambda$ and  no viscosity supersolutions for $\lambda < \bar\lambda$. It is also called ergodic constant or additive eigenvalue. The existence of such number is usually proved by the limit
\begin{equation}
\label{limdelta1}
\lim_{\de\to 0+} \de w_\de(x) = \bar\lambda 
\end{equation}
where 
\begin{equation}
\label{eqdelta1}
\de w_\de+ H(x, Dw_\de) =0 \quad\text{ in } \R^n , 
\end{equation}
the so-called vanishing discount problem. An alternative construction is the ergodic-type limit
\begin{equation}
\label{limt1}
\lim_{t\to +\infty} \frac{w(x,t)}t = \bar \lambda  .
\end{equation}
where
\begin{equation} 
\label{eq:2ter} 
\partial_t w +H(x, D_xw) = 0 ,\quad \text{in } 
\R^n\times \mathopen( 0,+\infty \mathclose) , \quad w(\cdot,0)= w_o\in C(\T^n) .
\end{equation} 
For Hamiltonians coming from calculus of variation or optimal control these results describe the asymptotics of  value functions and are related to ergodic control, as it was first observed by Lasry in 1975 \cite{Las75}, see \cite{BCD, Ari98, AL98}, \cite{ABdg, ABmem} for differential games, and the references therein.
The validity of one of the two limits \eqref{limdelta1} and \eqref{limt1} implies the other and it is called {\em ergodicity} of $H$.
Moreover, if for $P\in \R^n$ one considers the critical equation $\lambda+ H(x, Du+P )=0$, then the critical value $\bar \lambda = \bar H(P)$ is the effective Hamiltonian arising in periodic homogenization of Hamilton-Jacobi equations. This was the original motivation 
for studying  the critical value in the pioneering paper by Lions, Papanicolaou and Varadhan \cite{LPV}. 
In classical mechanics a smooth solution of the critical equation is a generating function of canonical transformations for integrable systems; in the literature on Hamiltonian systems $-\bar \lambda$ is called Ma\~n\'e critical value and $\bar H(P)$ is related to Mather's $\al$-function, see the surveys \cite{evans2004survey, fathi2008, Fa2014, sorr}.

The existence of the critical value via the limits \eqref{limdelta1} and \eqref{limt1} does not require the coercivity of $H$ but only a uniform bound on the oscillation of the solution $w_\de$ of \eqref{eqdelta1}
\begin{equation}
\label{bounddelta}
w_\de(x) -  w_\de(z) \leq C \quad
\forall\, x,z\in \T^n , \forall \, \de>0 .
\end{equation}
For Bellman Hamiltonians this follows from a
 weak controllability property of the underlying control system, called {\em bounded-time controllability} in \cite{ABmem}, see also \cite{Ari98, Gru98, AG00}
 . It also extends to 2nd order degenerate elliptic Hamiltonians arising in stochastic control \cite {AL98} and to the non-convex Isaacs Hamiltonians of differential games \cite{ABmem}. Barles \cite{barles07} proved the ergodicity of some nonconvex $H$ missing coercivity in one of the space dimensions by proving \eqref{bounddelta} with PDE methods. Cardaliaguet and Mendico \cite{CarM} proved \eqref{limt1} for calculus of variations problems where the control is the acceleration, instead of the velocity, so the Hamiltonian is linear in a half of components of $p$, and they apply this to Mean Field Games.

{\bf 2.} A stronger property of the critical value $\bar \la$ is that it is the unique $\la$ such that \eqref{eq:crit1} has a viscosity solution.  A function $\chi\in C(\T^n)$ solving
\begin{equation}
\label{criteq1}
\bar\lambda+ H(x, D\chi
 )=0 \quad\text{ in } \R^n 
\end{equation}
is called a {\em critical solution} of $H$ and a {\em corrector} in the jargon of homogenization. The existence of a Lipschitz critical solution was proved by Lions, Papanicolaou and Varadhan \cite{LPV} for coercive Hamiltonians on the torus, and by Fathi \cite{fathi1997} on compact manifolds without boundary for smooth $H$ convex conjugate of a Tonelli Lagrangian. This result is often called {\em weak KAM theorem}. In these papers the coercivity of the Hamiltonian gives the Lipschitz continuity of $\chi$, but it is not needed for a merely continuous solution if the Hamiltonian is of Bellman type. In fact, it was observed in \cite{ABmem} that for {\em small-time controllable} systems, i.e., such that for any $x, z\in\T^n$ there is a control driving the system from $x$ to $z$ in a time smaller than $\omega(|x-z|)$ with $\omega$ a modulus, there is a critical solution $\chi$ with modulus of continuity $C\omega$. This includes symmetric control systems whose vector fields 
are Lie-bracket-generating 
of order $k$, and then  $\omega(r)=Cr^{1/k}$ and $\chi$ is H\"older continuous (Example 2, p. 43 of \cite{ABmem}). In this case the Hamiltonian is \eqref{Ham1} with $b\equiv 0$, so it is bounded from below and coercive in the directions of the fields, it is sometimes called {\em pseudo-coercive}. See also \cite{BW, Stro} for similar results in Carnot groups, 
  \cite{Gom07} and \cite{CanMen23}, respectively, for vakonomic mechanics and for
   sub-Riemannian structures on the torus.

In \cite{LPV} a critical solution $\chi$ is constructed as the uniform limit of $w_\de(x) - \min_{\T^n} w_\de$ as $\de\to 0+$. The 
 two uniform limits
\begin{equation}
\label{twolim}
\lim_{t\to +\infty} \left(w(x,t) - \bar\lambda t \right)  = \chi(x)   \,, \qquad  \lim_{\de\to 0+}\left(w_\de(x) - \frac {\bar\lambda}\de \right) = \tilde\chi(x) \,,
\end{equation}
where $\tilde \chi$ is also critical, improve strictly \eqref{limdelta1} and \eqref{limt1} by giving an explicit connection between $\bar \la$ and critical solutions. They are both the subject of 
a large literature that we do not try to review here. 
Note that the former says that in the long time 
any solution of the evolutive HJ equation \eqref{eq:2ter}  
tends to the profile $\chi$ moving up or down with velocity $\bar \la$. Some relevant references about it are Fathi \cite{fathi1997, fathi1998},  Roquejoffre \cite{Roc},  and Barles and Souganidis \cite{BS} in the coercive case under different conditions, and Cannarsa and Mendico \cite{CanMen23} in the pseudo-coercive case. For the latter limit in \eqref{twolim} the main result was proved in \cite{DFIZ}.

{\bf 3.} Weak KAM theory also studies in depth the properties of the {\em Lax-Oleinik semigroup} on continuous functions $\phi$
\begin{equation}
\label{laxOl}
T_t \phi(x):= \inf_{z\in \T^n}\left\{ h_t(z,x)+\phi(z)\right\} 
\end{equation}
where $h_t
$ is the {\em minimal action} of the Lagrangian $L$ conjugate to $H$
\[
h_t(z,x) :=\inf
\left\{ \int_0^t  L(y(s), \dot y(s))\,ds \,:\, y(0)=z \,, y(t)=x \right\} .
\]
Another characterization of the critical value  is the existence of a common fixed point $\chi$ for $T_t - \bar\la t$ for all $t>0$, see Fathi \cite{
fathi1998, fathi2008, Fa2014} for the classical Tonelli setting, Cannarsa and Mendico \cite{CanMen23} for the sub-Riemannian case, and the references therein. The function $\chi$ is called {\em weak KAM solution} and it is also a critical solution. A crucial step for these results, as well as for the first limit in \eqref{twolim},  is an equicontinuity estimate on the functions $h_t$.

{\bf 4.} Homogenization problems concern HJ equations with highly oscillating Hamiltonian 
\begin{equation} 
\label{homo1} 
\partial_t v^\ep +H\left(z, \frac z\ep, D_zv^\ep \right) = 0 \,, \quad 
(z, t)\in \R^n\times \mathopen( 0, T \mathclose) , 
\end{equation} 
aiming at finding a limit to the solutions $v^\ep$ as $\ep\to 0+$. Lions, Papanicolaou and Varadhan \cite{LPV} showed, for coercive Hamiltonians 
$\Z^n$-periodic in the second entry, that the limit of $v^\ep$ satisfies the HJ equation 
\begin{equation} 
\label{eq;eff}
\partial_t v +\bar H(z,D_zv) = 0 \;\text{ in } \R^n\times \mathopen( 0,T \mathclose) ,
\end{equation} 
where the {\em effective Hamiltonian} $\bar H(z, P)$ is the critical value associated to the Hamiltonian $H(z, x, p + P)$ with $z$ and $P$  frozen. Among the many papers on the subject let us mention Evans \cite{E92}, who introduced the perturbed test function method that considerably simplified the proofs, and some papers for non-coercive Hamiltonians \cite{BW, barles07, Stro, ABmem, BT, CarNS, Sic, Vit}, see also the references therein. 

\medskip
Now we turn to Hamiltonians of the form \eqref{Ham1}. We will treat two cases: (i) the control set $A$ bounded,  so that $H$ has at most linear growth in $p$ at infinity, 
(ii) $A=\R^m$ and 
\begin{equation}
\label{Lag1}
L(x, a)=a^T G(x) a + l(x) ,
\end{equation} 
where $G(x)$ is positive definite and has an inverse with Lipschitz square root $\tau(x)$. Then the Hamiltonian is quadratic in $p$
\begin{equation}
\label{Hquad1}
H(x,p)= b(x)\cdot p + \frac12 |\sigma(x)p|^2 -l(x) , \quad \sigma(x):= (F\tau)^T(x) .
\end{equation}
In this Introduction we'll focus on case (ii). It was studied by Agrachev and Lee \cite{AgL} on general compact manifold without boundary $M$. They have two main results. The first assumes the family of vector fields $\{f^1,...,  f^m\}$ be 3-generating, namely, for each $x\in M$, $f^i(x)$ and the Lie brackets up to order 2, $[f^i,f^j](x)$ and $[f^k, [f^i,f^j]](x)$, span the tangent space to $M$ ($\R^m$ in our case). Then the minimal action $h_t$ is continuous, and in this case they can develop a weak KAM theory as in points 1, 2, and  3 above. The second main result of \cite{AgL} is, for any $k\geq 4$, an example of a  simple explicit family of two vector fields in two dimensions that is $k$-generating and such that $h_t$ is {\em discontinuous}. Hence some crucial estimates  of classical weak KAM theory fail and new methods appear to be necessary. 

The goal of the present paper is to go beyond this difficulty for 
general families of fields $f^i$ with the bracket generation condition, without restrictions on the order. 
We give results on all asymptotic properties described above, sometimes in a weaker form within the theory of discontinuous viscosity solutions.
To our knowledge this is the first time such theory is used for weak KAM problems.

In  Section \ref{existence} we prove the existence of the critical value $\bar\la$ and its approximation by the ergodic limits  \eqref{limdelta1} and  \eqref{limt1}, via the bounded-time controllability of the associated control system 
\begin{equation}
\label{syst1} 
\dot{y}(t)= 
 -b(y(t)) - \sum_{i=0}^m \alpha_i(t) f^i(y(t)),\quad t>0,
\end{equation}
with controls $\alpha(\cdot)$ taking values in a sufficiently large set of $\R^m$. This is proved in Sect. \ref{btc} by geometric control methods, in particular using a result of Brunovsky and Lobry \cite{BL}. Other useful tools throughout the paper, presented in Sect. \ref{compar}, are comparison principles between viscosity sub- and supersolutions, that are not standard for the quadratic non-coercive Hamiltonian \eqref{Hquad1} 
and are adapted  from \cite{BP90, Bal, BCG}. 

Section \ref {weakK} is devoted to several results of weak KAM type.  In Section \ref{existence2} we show that 
the critical equation \eqref{criteq1} has an u.s.c. subsolution, a l.s.c. supersolution, and not-necessarily-continuous 
 solutions. For various definitions of possibly discontinuous viscosity solutions see Sect. \ref{existence2} and the survey in Chapt. V of \cite{BCD}.

Section \ref{asymptotic} provides some weakened
 versions of the limits \eqref{twolim} in terms of relaxed semilimits. The main result on the long-time behavior of solutions to the evolutive equation \eqref{eq:2ter}  is that
 \[
 v_{\#}(x) := \liminf_{t\to +\infty, y\to x} \left(w(y, t) - \bar\lambda t \right)
 \]
 is a lsc-solution, or bilateral supersolution, of \eqref{criteq1}, namely,
 \[
 \bar\la + H(x, p) = 0 \quad \forall\, p\in D^-v_{\#}(x) \quad \forall\, x\in\R^n \,,
 \]
 where $D^-$ denotes the usual sub-differential of a lsc function.
 
 In Section \ref{sec:fixed} we consider the Hamilton-Jacobi semigroup
 \[
T_t \phi(x):= \inf_{\alpha\in\mathcal{A}}\left\{ \int_0^t  L(y(s), \alpha(s))\,ds+\phi(y(t)) : y(\cdot) \text{ solves } 
 \eqref{syst1}  , 
 y(0)=x \right\} 
\]
which is a forward variant of the Lax-Oleinik semigroup \eqref{laxOl}. We study the existence of weak KAM solutions, i.e., fixed points of $T_t - \bar\la t$, and their connection to the critical equation. The main result is the existence of a function $\chi$ such that
  \begin{equation*}  
\chi(x) + \bar\lambda t = T_t \chi (x) 
 \quad \forall x\in \R^n, \forall  t\geq 0 ,
\end{equation*}
the usc envelope $\chi^*$ is a subsolution of  \eqref{criteq1}, and the lsc envelope $\chi_*$ is a  bilateral supersolution of  \eqref{criteq1}.
 
 Section \ref{homo} is an application of the results on the critical value in Section \ref{existence} to the periodic homogenization of stationary, discounted HJ equations as well as of evolutive equations. We use the methods of \cite{E92, ABsiam, ABarma, ABmem} to prove the uniform convergence of solutions of \eqref{homo1} to the solution of \eqref{eq;eff} satisfying suitable initial data, where $\bar H$ is the critical value of an auxiliary Hamiltonian with parameters.
 
 In this paper we dot treat the Aubry-Mather theory that completes the weak KAM theory in the classical cases \cite{fathi2008, Fa2014} and recently extended to the sub-Riemannian setting in \cite{CanMen23}. It is a possible natural continuation of the present research. Another one is dropping the compactness of the state space and considering problems in all $\R^n$ without periodicity. Several classical results were extended in that direction, see, e.g., the recent papers \cite{ishiiSic,  CanMen22, BK} and the references therein.
 
 \smallskip {\em Acknowledgements}. This research originates 
 from the joint work with Sarah Balbinot for her Master thesis \cite{Bal},  some results of Sections \ref{Prelim}, \ref{existence}, and \ref{homo} appeared in  it.  
\section{Basic assumptions and 
preliminary results} 
\label{Prelim} 
We consider Hamiltonians of the form
\begin{equation}
\label{Ham}
H(x,p)= b(x)\cdot p + \sup_{a\in A} \left\{ \sum_{i=1}^m a_i f^i(x) \cdot p - L(x, a) \right\},
\end{equation}
where $A\subseteq\R^m$, $a=(a_1,...,a_m)$, under the standing assumptions that $b, f^i : \R^n \to\R^n$, $i=1,...,m$, are Lipschitz and $\mathbb{Z}^n$-periodic, i.e., 
\[ 
b(x+k)=b(x),\quad f^i(x+k)=f^i(x),\quad \forall k\in\mathbb{Z}^n ,
\]
and $L : \R^n\times \bar A \to \R$ is continuous and $\mathbb{Z}^n$-periodic in $x$. We will denote with $\mathcal{M}_{l,k}$ the $l\times k$ real matrices and with $M^T$  the transpose of a matrix $M$. We also call $F(x)\in \mathcal{M}_{n,m}$ the matrix whose columns are the vector fields $f^i$, i.e.,
\begin{equation}
\label{F}
F(x)_{ji}=f^i_j(x) ,\qquad j=1,..., n, \, i=1,..., m ,
\end{equation}
 so that the summation in the definition of $H$ can be written $p\cdot F(x)a$.
We will treat two cases: 
 (i) $A$ bounded, 
 (ii) $A=\R^m$ and 
\begin{equation}
\label{Lag}
L(x, a)=a^T G(x) a + l(x) ,
\end{equation} 
where $G(x)$ is a positive definite $m\times m$ symmetric matrix, 
and its inverse matrix 
 has Lipschitz square root, $\mathbb{Z}^n$-periodic,
\begin{equation}
\label{Ginv}
G^{-1}(x) = \tau(x) \tau^T(x) , \qquad \tau \colon \R^n\to \mathcal{M}_{m,m} \quad\text{ Lipschitz}.
\end{equation}
Since $G^{-1}$ is positive definite, the sup in the definition of $H$ is a max, attained at $a=-G^{-1}F^Tp$, and then
\[
H(x,p)= b(x)\cdot p + \frac12 p^TFG^{-1}F^Tp -l(x) 
\]
which can rewritten as 
\begin{equation}
\label{Hquad}
H(x,p)= b(x)\cdot p + \frac12 |\sigma(x)p|^2 -l(x) , 
\end{equation}
where 
\begin{equation}
\label{sig}
\sigma(x):= (F\tau)^T(x)
\end{equation}
 is Lipschitz, $l$ is continuous and  both are periodic.
The Hamiltonian \eqref{Ham} is associated to control problems with running cost, or Lagrangian, $L(x, a)$, and trajectories driven by the affine system
\begin{equation}
\label{syst} 
\dot{y}(t)=  -b(y(t)) - F(y(t)) \al(t) = -b(y(t)) - \sum_{i=0}^m \alpha_i(t) f^i(y(t)),\quad t>0,
\end{equation}
with control functions $\alpha(\cdot)\in \mathcal{A}:= L^\infty( \mathopen[ 0,T\mathclose], A)$. 
 
 In the rest of the paper, when we refer to "the case $A=\R^m$" we mean that also \eqref{Lag}, \eqref{Ginv}, \eqref{Hquad}, and \eqref{sig} hold true.

\smallskip
\noindent
{\bf Notations}. We say that a function $\phi : X\to\R$, $X\subseteq\R^k$, is in $B(X)$ or, respectively, $BC(X)$, $BLSC(X)$, $BUSC(X)$, $BUC(X)$, if is bounded or, respectively, bounded and continuous, bounded and lower semicontinuous (briefly, lsc), bounded and upper semicontinuous  (briefly, usc), bounded and uniformly continuous.

\subsection{Comparison principles}
\label{compar}
We say that a Hamiltonian $H$ {\em satisfies the comparison principle} if, for all $\delta>0$, any subsolution $u\in BUSC(\R^n)$ and supersolution $v\in BUSC(\R^n)$ of
\[
\delta w +H(x, Dw) = 0 \qquad \text{in } \R^n 
\]
verify $u \leq v$ in $\R^n $.

This property is well known for  $H$ defined by \eqref{Ham} in the  case of linear growth, i.e., $A$ bounded, see, e.g., \cite{barles, BCD}. For the quadratic case it requires some additional tricks introduced in \cite{BP90, barles}.
\begin{thm}
\label{teo:comp}
Assume $H$ is given by \eqref{Hquad} with $b, \sigma$ Lipschitz and bounded and $l\in BC(\R^n)$. Then the comparison principle holds.
\end{thm}
\begin{proof} It is obtained by adapting with standard methods the proof of Thm. 11 in \cite{BP90} or the proof of Thm. 3.5 in \cite{BCG}, see also Thm. 1.1.5 in \cite{Bal}.
\end{proof}
Next we state the comparison principle for the evolutive HJ equation
\begin{equation}
\label{eq:2} 
\partial_t w +H(x, Dw) = 0 ,\quad (x,t)\in\R^n\times \mathopen( 0,T\mathclose).
\end{equation}
\begin{thm}
\label{teo:confronto}
Assume $H$ is given either by \eqref{Ham} with $A$ bounded or by \eqref{Hquad} with with $b, \sigma$ Lipschitz and bounded and $l\in BC(\R^n)$.
Let $u\in BUSC(\R^n\times \mathopen[ 0,T\mathclose]),v\in BLSC(\R^n\times \mathopen[ 0,T\mathclose])$ be, respectively, a sub- and a supersolution of~\eqref{eq:2}. Then
\[
\sup_{\R^n\times \mathopen[ 0,T\mathclose]} ( u-v)\leq \sup_{\R^n}( u-v)(\cdot, 0).
\]
\end{thm}
\begin{proof} The case of $A$ bounded is well-known \cite{barles, BCD}. For  the quadratic Hamiltonian \eqref{Hquad} with $b\equiv 0$ the result follows immediately from  Thm. 3.5 in \cite{BCG}.  The proof can be adapted to the case $b\neq 0$ 
 by standard methods, see Thm. 1.1.2 in \cite{Bal} for the details.
\end{proof}

\subsection{Bounded-time controllability of affine control systems}
\label{btc}
We consider the system \eqref{syst}
and denote with $y_x(\cdot, \al)$ the trajectory starting at $x$, i.e., with $y_x(0, \al)=x$. We will also denote with $B(x, R)$ the open ball in $\R^n$ with center $x$ and radius $R$, and
\[
B_m(R):=\{ a\in\R^m : |a| 
< R \} .
\]

Define the time taken by the system to join two points as 
\[
t^\sharp(x_1, x_2):=\inf\{ t\geq 0\,|\, \exists \alpha\in\mathcal{A}\text{ such that } y_{x_1}(t, \alpha) - x_2 \in \Z^n 
\},
\]
which is $+\infty$ if there is no trajectory joining the points. If $t^\sharp$ is finite for all pairs of points the system is called completely controllable in \cite{CK}.

Our main assumption is the {\em strong Lie bracket generation condition}: 
\begin{equation}
\label{SBG}
\tag{SBG}
\text{the vector fields } f_1,...,f_m \text{ and their commutators of any order span $\R^m$ at each point of  } \R^m .
\end{equation}
The next complete controllability result is obtained from Prop. III.8 in the paper \cite{BL} by Brunovski and Lobry in the special case of state space the flat torus $\T^n=\R^n / \Z^n$, which is compact with the norm $|x_1 - x_2|_T:= \inf_{k\in \Z^n} |x_1 - x_2+k|$.

\begin{prop}
\label{teo:BL} 
(Brunovski and Lobry)  Assume (SBG). Then there exists $K>0$ such that, if $A\supseteq B_m(K)$, then $t^\sharp(x_1, x_2)<+\infty$ for all $x_1, x_2$.
\end{prop}

We will also need an estimate on $t^\sharp$. Following the terminology of \cite{ABmem} we say that the system \eqref{syst} is {\em bounded-time controllable} on $\T^n$ if 
\begin{equation}
\label{BTC} 
\tag{BTC} 
\text{there exists  }\, S>0 \quad\text{such that} \quad t^\sharp(x_1, x_2)\leq S \quad\text{for all } x_1, x_2\in \R^n.
\end{equation}
This property is also called {\em uniform exact controllability} in \cite {Ari98} and {\em total controllability} in \cite{AG00}.

The next result follows from Krener's theorem stating the local accessibility of \eqref{syst} under a weak bracket generation condition (implied by \eqref{SBG}), 
see Thm. A.4.4 in  \cite{CK}.
\begin{prop}
\label{teo:Krener} 
Assume (SBG). Then for all $\bar x$, the set
\[
\mathopen\{ {x}\,|\, \exists \,0< t\leq 1, \alpha\in\mathcal{A} \text{ such that }   y_x(t, \alpha) =\bar x 
\mathclose\}
\]
has non-empty interior.
\end{prop}

The next lemma combines the two previous propositions and it will be crucial for proving the bounded-time controllability of \eqref{syst}.
\begin{lemma}
\label{halfsbg}
Assume (SBG) and $A\supseteq B_m(K)$. Then for any $\bar x$ there exists 
$C(\bar x)$ such that
\begin{equation}
\label{halfest}
t^\sharp(x, \bar x) \leq C(\bar x) \qquad \forall\, x .
\end{equation}
\end{lemma}
\begin{proof}
By Thm. \ref{teo:Krener} there exist $x_1$ and $r_1>0$ such that for all $z\in B(x_1, r_1)$ there is a control $\al_z$ and $t\in(0,1]$ satisfying
\begin{equation}
\label{uno}
y_z(t, \al_z) = \bar x .
\end{equation}
On the other hand, for any $x$, Thm. \ref{teo:BL} provides a control $\al$ 
 and $T>0$ such that  
\[
y_x(T, \al) =  x_1 .
\]
Note that  $\al$ and $T$ depend on $x$ and $x_1$, and $x_1$ depends only on $\bar x$, so we will write $\al=\al_{x, \bar x}$ and $T=T(x, \bar x)$.
Now we use the continuous dependence on the initial position of the solution of \eqref{syst} with fixed control $\al_{x, \bar x}$ and get $r_2=r_2(x, \bar x)$ such that 
\[
y_w(T, \al_{x, \bar x}) \in B(x_1, r_1) \qquad \forall \, w\in B(x, r_2) .
\]
Next we take the concatenation $\bar \al$ of the controls $\al_{x, \bar x}$ until time $T(x, \bar x)$ and $\al_z(\cdot-T)$, with $z=y_w(T, \al_{x, \bar x})$, and obtain
\[
y_w(T+t, \bar \al)= y_z(t, \al_z) = \bar x , \qquad \forall\, w\in B(x, r_2) ,
\]
by \eqref{uno} and using that $z\in B(x_1, r_1)$.

Now we use the compactness of $\T^n$ to extract a finite number of point $x_j$, $j=1,...,,N$ such that the balls $B(x_j, r_2(x_j, \bar x))$ cover the unit cube. Then, for any $y\in\T^n$ there is $k$ such that $y\in B(x_k, r_2(x_k, \bar x))$ and so 
\[
t^\sharp(y, \bar x) \leq T(x_k, \bar x) +t \leq \max_j T(x_j, \bar x) +1 =: C(\bar x) .
\]
\end{proof}

The next is the main result of this section.
\begin{thm}
\label{teo:btc} 
Assume \eqref{SBG} and $A\supseteq B_m(K)$. Then the system \eqref{syst} has the bounded-time controllability property \eqref{BTC}.
\end{thm}
\begin{proof}
Consider the system
\begin{equation*}
\dot{y}(t)=  b
(y(t)) + \sum_{i=0}^m \alpha_i(t) f_i(y(t)),\quad t>0,
\end{equation*}
whose trajectories are the trajectories of \eqref{syst} run backward. This system also satisfies the bracket condition \eqref{SBG}. Then we can apply Lemma \ref{halfsbg} to it and get, for any $\bar x$,
\[
t^\sharp(\bar x, x) \leq C(\bar x) \qquad \forall\, x .
\]
Combining this estimate with \eqref{halfest} we obtain, for any $x_1, x_2$
\[
t^\sharp(x_1, x_2) \leq t^\sharp(x_1, 0) + t^\sharp(0, x_2) \leq 2 C(0) ,
\]
which proves the claim.
\end{proof}

\section{The critical value: existence and approximation}
\label{existence}

The critical equation for a Hamiltonian $H$ is
\begin{equation}
\label{eq:cella}
\lambda+ H(x, Du
 )=0 
\end{equation}
where the unknowns are $\lambda\in \R$, called critical value, and a periodic function $u$, 
called critical solution, or corrector.   
The candidate critical values are
\[ 
\lambda_1:= \sup\{ \lambda : 
 \exists \text{ a subsolution in } BUSC(\R^n) 
\text{ of}~\eqref{eq:cella}\},
\] 
\[ 
\lambda_2 := \inf\{ \lambda  :  \exists \text{ a supersolution in } BLSC(\R^n) 
\text{ of}~\eqref{eq:cella}\}.
\]
The next lemma is well-known, see \cite{LPV, AL98, ABarma, ABmem}.
\begin{lemma}
If $H$ satisfies the comparison principle, then $\lambda_1\leq\lambda_2$.
\end{lemma}
\label{l1l2}
The next  is the main result of the present Section \ref{existence}. It identifies the critical value and approximates it by solutions of
 the stationary discounted HJ equation 
\begin{equation}
\label{eqdelta}
\de w_\de+ H(x, Dw_\de) =0 \quad\text{ in } \R^n . 
\end{equation}
In the terminology of \cite{ABarma, ABmem} it says that the Hamiltonian $H$ is {\em ergodic}. It gives enough informations for the applications to homogenization of 
Section \ref{homo}. However it can be made more precise and more similar to weak KAM theorems, see the next Section \ref{weakK}.

The statement and proof involve the auxiliary "truncated" Hamiltonian
\[
H_K(x,p) := b(x)\cdot p + \sup_{a\in B_m(K)} \left\{ p\cdot F(x) a 
 - L(x, a) \right\}
\]
where $K$ is the constant appearing in Thm. \ref{teo:BL} and the data $b, F, $ and $L$ satisfy the assumptions of Section \ref{Prelim}.
\begin{prop}
\label{conv1}
Let $H$ be any continuous Hamiltonian satisfying the comparison principle, $|H(x,0)|\leq C$ for all $x$, and $H(x,p)\geq H_K(x,p)$ for all $x, p$. Assume the vector fields $f^1,..., f^m$ satisfy \eqref{SBG} and, for each $\de>0$, $w_\de$ solves \eqref{eqdelta}. Then $\lambda_1=\lambda_2=:\bar\lambda$ and
\[
\lim_{\de\to 0} \de w_\de(x) = \bar\lambda \qquad\text{uniformly in } \R^n .
\]
\end{prop}
\begin{proof}
Since $w_\de$ are subsolutions of $w_\de+ H_K(x, Dw_\de) =0$, they 
 satisfy a suboptimality principle \cite{BCD}(Thm. III.2.32), namely, for all $\tau>0$ 
 \[
 w_\de(x) \leq \inf_{\alpha\in\mathcal{A}_K}\left\{ \int_0^\tau e^{-t\delta} L(y_x(t, \al), \alpha(t))\,dt+e^{-\tau \delta}w_\delta(y_x(\tau, \al))\right\} ,
 \]
 where $\mathcal{A}_K$ denotes the measurable functions of time taking values in $B_m(K)$.
We fix $x, \bar x$ 
  and use the property \eqref{BTC} from Thm. \ref{teo:btc} to find $\alpha\in\mathcal{A}_K$ and 
  $\bar t\leq 2 t^\sharp(x, \bar x)
  $ such that $y_x(\bar t, \al)-\bar x\in \Z^n$. Call $C_L:= \max_{\R^n\times B_m(K)} |L|$. Then, taking $\tau=\bar t$,
\[
w_\delta(x)-e^{-\bar t \delta}w_\delta(\bar x)\leq\frac{ C_L}{\delta}(1-e^{-\delta 2S
)}) .
\]
The condition $|H(x,0)|\leq C$ implies that $-C/\de$ and $C/\de$ are, respectively, a sub and a supersolution of  \eqref{eqdelta}. Then the comparison principle gives 
\[
-C/\de \leq w_\de(x) \leq C/\de \qquad \forall \, x 
\]
and
\[
\de w_\delta(x) -\de w_\delta(\bar x) \leq 
(C_L + C) (1-e^{-\delta 2S 
}) \qquad \forall\, x, \bar x.
\]
By exchanging the roles of $x$ and $\bar x$ we obtain
\begin{equation}
\label{bigO}
|\de w_\delta(x) - \de w_\delta(\bar x)| \leq  O(\de)  \qquad \forall\, x, \bar x.
\end{equation}
Now we fix $\bar x$ and extract a sequence $\de_k\to 0$ such that $\de_k w_{\delta_k}(\bar x) \to \mu \in \R$. Then \eqref{bigO} implies the uniform convergence of $\de_k w_{\delta_k}$ to $\mu$. 

In order to show that the limit is unique we take $\lambda < \mu$ and an index $k$ such that $\lambda<\de_k w_{\delta_k}(x)$ for all $x$. Then $w_{\delta_k}$ is a subsolution of \eqref{eq:cella} and so $\lambda \leq \lambda_1$. This gives $\mu \leq \lambda_1$. A symmetric argument with supersolutions gives $\mu \geq \lambda_2$, which combined with Lemma \ref{l1l2} completes the proof.
\end{proof}
\begin{cor}
\label{conv2}
For $H$ defined by \eqref{Ham} and under the assumptions of Section \ref{Prelim}, suppose either that $A\supseteq B_m(K)$ is bounded, or that $A=\R^m$ and
$$
L(x, a)=a^T G(x) a + l(x)$$
 with $G$ and $l$ as in Section \ref{Prelim}. Suppose the vector fields $f^1,..., f^m$ satisfy \eqref{SBG} and $w_\de$ solves \eqref{eqdelta}. Then the conclusions of Prop. \ref{conv1} hold true.
 \end{cor}
\begin{proof} The assumptions of Prop. \ref{conv1} are easily checked if $A$ is bounded. In the case of quadratic $L$ the comparison principle is Thm. \ref{teo:confronto} and 
$H(x,0) = -l(x)$
which is bounded, so again Prop. \ref{conv1} gives the conclusion.
\end{proof}
\begin{rem}
Note that Proposition \ref{conv1} 
applies also to Hamiltonians not necessarily of the form \eqref{Ham} and possibly also non-convex in $p$.
\end{rem}
Next we approximate the critical value by the 
long time limit of the solution of the evolutive H-J equation \eqref{eq:2}. This is a result of Abelian-Tauberian type similar to Thm. 5 in \cite{Ari98} and Thm. 3 in \cite{ABarma}.
\begin{prop}
\label{convtime}
Assume the same conditions of Cor. \ref{conv2} and $w(x,t)$ solution of the Cauchy problem 
\begin{equation} 
\label{eq:2bis} 
\partial_t w +H(x, D_xw) = 0 ,\quad (x,t)\in\R^n\times \mathopen( 0,+\infty \mathclose) , \quad w(\cdot,0)= w_o ,
\end{equation} 
with $w_o\in C(\R^n)$ and $\Z^n$-periodic. Then 
\[
\lim_{t\to +\infty} \frac{w(x,t)}t = \bar \lambda  = \lambda_1 = \lambda_2
 \quad\text{uniformly in } x .
\]
\end{prop}
\begin{proof}
Let $\mu\in \R$ be such that there is a subsolution $v\in BUSC(\R^n)$ of $\mu +H(x, Dv)\leq 0$. By subtracting a constant it is not restrictive to assume $v\leq w_o$. Then $v+\mu t$ is a subsolution of the Cauchy problem \ref{eq:2bis} . The comparison principle for such problem is standard for $A$ bounded and it is Thm. \ref{teo:confronto} for the quadratic case. Then $v(x)+\mu t\leq w(x, t)$ for all $x$. We divide by $t$, send $t$ to $+\infty$ and take the sup over $\mu$ to get
\[
\liminf_{t\to +\infty} \min_x \frac{w(x, t)}t \geq \lambda_1 . 
\]
In a similar way we get 
\[
\limsup_{t\to +\infty} \max_x \frac{w(x, t)}t \leq \lambda_2 
\]
and the conclusion follows from Prop. \ref{conv1}.
\end{proof}
\section{Weak KAM theorems}
\label{weakK}
In all this section we consider the critical equation \eqref{eq:cella} for the Hamiltonian defined by \eqref{Ham1} (equivalently, by  \eqref{Ham}), under the assumptions of Corollary \ref{conv2}, namely,
\begin{itemize}
\item $b, f^i$ are Lipschitz and $\Z^n$-periodic, the columns $f^i$ of $F$ verify the bracket condition \eqref{SBG};
\item 
 either 
(i) $A$ is bounded and $A\supseteq B_m(K)$, 
\item or (ii) $A=\R^m$, $L$ is given by \eqref{Lag} with $G$ 
satifying \eqref{Ginv}. 
\end{itemize}
In the 
 case (ii) the critical equation \eqref{eq:cella} can be written as
\begin{equation} 
\label{eq:crit}
\lambda + b(x)\cdot Du + \frac12 |\sigma(x)Du|^2 = l(x) \quad \text{in } \R^n ,
\end{equation} 
 with $\sigma$ given by \eqref{sig}.
 
\subsection{Existence of a critical solution}
\label{existence2}
We recall two definition of viscosity solution valid for not necessarily continuous functions, the first due to Ishii \cite{ishiiPerron} and the second to Barron and Jensen \cite{BJ}, see also \cite{Fran}. A survey is in \cite{BCD}, Sections V.2 and V.5, respectively. 

{\bf Definitions.} A locally bounded function $u$ is a {\em non-continuous viscosity solution} of a PDE in an open set if its upper semicontinuous envelope $u^*(x):=\limsup_{y\to x} u(y)$ is a viscosity subsolution and its lower semicontinuous envelope $u_*(x):=\liminf_{y\to x} u(y)$ is a viscosity supersolution.

A function $u\in LSC(\Omega)$ is a {\em lower semicontinuous viscosity solution}, briefly, {\em lsc-solution}, or {\em bilateral supersolution} of $F(x,u,Du)=0$ in $\Omega$ open set, if $F(x, u(x), p)=0$ for all $p\in D^-u(x)$ and $x\in\Omega$, where $D^-u(x)$ is the subdifferential of $u$ at $x$. In other words, $u$ is a viscosity supersolution of both $F=0$ and  $-F(x,u,Du)=0$.
\begin{thm}
\label{kam}
Under the standing assumptions of this section
\begin{equation} 
\label{maxmin}
\begin{split}
\bar \lambda &= \max\{ \lambda : 
 \exists \text{ a subsolution in } BUSC(\R^n) 
\text{ of}~\eqref{eq:cella}\},\\
&= \min\{ \lambda  :  \exists \text{ a supersolution in } BLSC(\R^n) 
\text{ of}~\eqref{eq:cella}\} ,
\end{split}
\end{equation}
and so $\lambda=\bar \lambda$ it the unique constant such that the critical equation \eqref{eq:cella} has both a sub- and a supersolution.

Moreover, for $\lambda=\bar \lambda$ the equation has a non-continuous viscosity solution and a lsc-solution. 
\end{thm}
\begin{rem}
In Section \ref{sec:fixed} we prove that there exists also a function $\chi\in B(\R^n)$ such that $\chi$ is a non-continuous viscosity solution and $\chi_*$ is a lsc-solution of \eqref{eq:cella}. 
\end{rem}
\begin{proof}
Define $v_\de(x):= w_\de(x)-w_\de(0)$, where $w_\de$ solves \eqref{eqdelta}. By the estimate \eqref{bigO} we have $|v_\de(x)|\leq \bar C$ for all $x$. Moreover $v_\de$ solves
\begin{equation}
 \label{eq:v}
\de v_\de + \de w_\de(0) + H(x, Dv_\de) = 0 .
\end{equation}
Note that $\de v_\de\to 0$ as $\de\to 0$ and $\de w_\de(0)\to \bar \lambda$ by Cor. \ref{conv2}.
Consider the relaxed semilimits
\[
\bar v(x) := \limsup_{\de\to 0, y\to x} v_\de(y) ,\quad \underline v(x) := \liminf_{\de\to 0, y\to x} v_\de(y) \,.
\]
By the well-known stability properties of viscosity semi-solutions we get
\[
 \bar \lambda +H(x, D\bar v)\leq 0 , \quad \bar \lambda +H(x, D\underline v)\geq 0 ,
\]
and so the sup and inf in the definitions of $\lambda_1$ and $\lambda_2$ are attained, which proves \eqref{maxmin}.

The existence of a non-continuous solution is proved by the Perron-Ishii method \cite{ishiiPerron}. Observe that $|\bar v(x)|, |\underline v(x)|\leq \bar C$ for all $x$, and so $\bar v - 2 \bar C\leq \underline v$. Since $\bar v - 2 \bar C$ is a subsolutions of \eqref{eq:cella} and $\underline v$ is a supersolution, by Thm. V.2.14 in \cite{BCD} we get that
\[
u(x):= \sup \{ v(x) : \bar v - 2 \bar C\leq v \leq \underline v \,, \, v^* \text{ subsolution of }  \eqref{eq:cella} \}
\]
is a non-continuous solution of  \eqref{eq:cella}.

To prove the existence of a lsc-solution we claim that $v_\de$ is also a bilateral supersolution of \eqref{eq:v}. Then the semilimit $\underline v\in BLSC(\R^n)$ satisfies also
\[
-\bar \lambda - H(x, D\underline v)\geq 0 ,
\]
which is what we want. To prove the claim we recall that $v_\de$ is the value function of an optimal control problem
 \[
 v_\de(x) = \inf_{\alpha\in\mathcal{A}}
 \int_0^\infty e^{-t\delta} \left(L(y_x(t, \al), \alpha(t)) - \de w_\de(0)\right) \,dt \,.
  \]
Then $v_\de$ satisfies a backward dynamic programming principle, see Prop. III.2.25 in \cite {BCD}, which implies the claim, as in Cor. III.2.28 of \cite {BCD}. 
\end{proof}
\begin{rem}
\label{contchi}
If $b\equiv 0$ then there is a continuous solution of \eqref{eq:cella} with $\la=\bar\la$, provided that $A$ contains a neighbourhood of $0$, even smaller than $B_m(K)$. In fact, in this case the system \eqref{syst} is symmetric and by \eqref{SBG} it is small time 
controllable everywhere. This means that \eqref{BTC} is strengthened to the existence of a modulus $\omega$ such that
\[
t^\#(x_1, x_2)\leq \omega(|x_1-x_2|_T) \quad \forall \,x_1, x_2 \in \R^n .
\]
Then in the proof of Prop. \ref{conv1} the estimate \eqref{bigO} is improved to
\[
|\de w_\de(x) - \de w_\de(\bar x)| \leq O(\de) \omega(|x-\bar x|_T)  . 
\]
This implies that $w_\de$ are equicontinuous, and therefore so are $v_\de$ in the proof of Thm. \ref{kam}. Then by Ascoli-Arzela there is a sequence $v_{\de_j}$ converging uniformly to a continuous $v$, which solves the critical equation by the stability properties of viscosity solutions. In the case of $A$ bounded this was proved in Rmk. 6.2 and Example 2, p. 43 of \cite{ABmem}.
\end{rem}
\subsection{Asymptotic behaviour of solutions to HJ equations}
\label{asymptotic}
The proof of  Theorem \ref{kam} 
 constructs  critical sub- and supersolutions from the solution of the discounted equation \eqref{eqdelta}. The next results gives a more explicit connection among the critical value, critical sub- and supersolutions, and the approximating equations \eqref{eqdelta} as well as the evolutive HJ equation \eqref{eq:2bis}. 
We recall that in classical 
weak KAM theory for coercive Hamiltonians the limits \eqref{twolim} exist
uniformly, and they produce critical solutions, i.e., they solve the critical equation 
\begin{equation} 
\label{eq:crit5}
\bar\lambda + H(x, Du) = 0   \quad \text{in } \R^n ,
\end{equation} 
which for the case (ii), $A=\R^m$, becomes
\begin{equation} 
\label{eq:crit2}
\bar\lambda + b(x)\cdot Du + \frac12 |\sigma(x)Du|^2 = l(x) \quad \text{in } \R^n .
\end{equation} 
In the current context we do not expect such a result because continuous critical solutions may not exist. We can prove the following weaker result.
\begin{thm}
\label{kam2}
Under the assumptions of Thm. \ref{kam}, the semilimits of solutions of the discounted equation \eqref{eqdelta},
\begin{equation} 
\label{semilim1}
\bar w(x) := \limsup_{\de\to 0, y\to x} \left(w_\de(y) - \frac {\bar\lambda}\de \right)  \quad\text{and}\quad \underline w(x) := \liminf_{\de\to 0, y\to x}  \left(w_\de(y) - \frac {\bar\lambda}\de \right) ,
\end{equation}
are, respectively, a bounded sub- and supersolution of the critical equation \eqref{eq:crit5}. 

The same conclusion holds for
the semilimits of solutions of the evolutive  equation \eqref{eq:2bis},
\begin{equation} 
\label{semilim2}
v^{\#}(x) := \limsup_{t\to +\infty, y\to x} \left(w(y, t) - \bar\lambda t \right) \quad\text{and}\quad   v_{\#}(x) := \liminf_{t\to +\infty, y\to x} \left(w(y, t) - \bar\lambda t \right) .
\end{equation}

Moreover, $\underline w$ and  $v_{\#}$ are lsc-solutions of \eqref{eq:crit5}. 
\end{thm}
\begin{proof}
Define $u_\de (x):= w_\de(x) - \bar\lambda/\de$ and observe that it solves 
\begin{equation} 
\label{eq:v2}
\de u_\de + \bar\lambda + H(x, Du_\de) = 0 .
\end{equation} 
We know from Cor. \ref{conv2} that $\de u_\de\to 0$ as $\de\to 0$. We claim that the family $u_\de$ is equibounded.
Consider the sub- and supersolutions $\bar v$ and $\underline v$ found in Thm. \ref{kam} and a constant $C$ such that $\bar v-C\leq 0, \underline v + C\geq 0$. Then $\bar v-C$ and $\underline v + C$ are, respectively, a sub- and a supersolution of \eqref{eq:v2}. By the comparison principle Thm. \ref{teo:comp} we get
\[
\bar v-C \leq u_\de \leq \underline v + C ,
\]
and then, for some $\bar C$, $|u_\de(x)|\leq \bar C$ for all $x$ and $\de>0$. Therefore the semilimits $\bar w$ and $\underline w$ of $u_\de$ are bounded and so they are, respectively, a sub- and a supersolution of \eqref{eq:crit2}, by applying the stability property of viscosity solutions to equation \eqref{eq:v2}.

Moreover, as in the previous proof, $u_\de$ is the value function of an infinite horizon optimal control problem, so it satisfies a backward dynamic principle and it is a bilateral supersolution of \eqref{eq:v2}. Then the semilimit $\bar w$ solves also $-\bar\lambda -H(x, D\underline w)\geq 0$.

For the second statement define $z_\eta(x,t):= w(x , \eta t) - \bar\lambda \eta t$, where $w$ solves \eqref{eq:2bis}, and note that it solves  
\begin{equation} 
\label{eq:z}
\frac1\eta\partial_tz_\eta + \bar\lambda + H(x, Dz_\eta) = 0 ,\quad (x,t)\in\R^n\times \mathopen( 0,+\infty \mathclose) . 
\end{equation} 
To prove the equiboundedness of $z_\eta$, consider again  $\bar v$ and $\underline v$ as above and pick a constant $C$ such that $\bar v-C\leq w_o$ and $\underline v + C \geq w_o$. 
Then $\bar v-C + \bar\lambda t$ and $\underline v + C + \bar\lambda t$ are, respectively, a sub- and a supersolution of the Cauchy problem \eqref{eq:2bis}, so the comparison principle Thm. \ref{teo:confronto} gives
\[
\bar v(x)-C + \bar\lambda t \leq w (x,t) \leq \underline v(x) + C + \bar\lambda t .
\]
Therefore there is $\bar C$ such that $|z_\eta(x,t)|\leq \bar C$  for all $(x,t)\in \R^n\times [0, +\infty)$ and $\eta>0$. This implies the boundedness of the semilimits 
\[
\bar z(x, t) := \limsup_{\eta\to +\infty, y\to x, s\to t} z_\eta(y, s) \quad\text{and}\quad   \underline z(x,t ) := \liminf_{\eta\to +\infty, y\to x, s\to t} z_\eta(y, s) .
\] 
Now the stability property of viscosity solutions for equation \eqref{eq:z} as $\eta\to +\infty$  gives that, for $t>0$, $\bar z(x, t)$ and $ \underline z(x,t)$ are a sub- and a supersolution of \eqref{eq:crit5}. Finally we observe that, for all $t>0$, $\bar z(x, t) =v^{\#}(x)$ and $\underline z(x, t) =v_{\#}(x)$, which gives the conclusion.

It remains to prove that $v_{\#}$ satisfies also $-\bar\lambda -H(x, Dv_{\#})\geq 0$. Since $z_\eta$ is the value function of a finite horizon control problem, it satisfies a backward DPP as in Prop. III.3.20 of \cite{BCD}. Then it is a bilateral supersolution of \eqref{eq:z}, by Thm. III.3.22 of \cite{BCD}. Therefore its lower semilimit  $v_{\#}$ verifies the desired inequality.
\end{proof}
\begin{rem}
Note that our result on the long-time limits \eqref{semilim2} is weaker that Theorem 9.14 of \cite{CanMen23} for the pseudo-coercive case without drift, based on equicontinuity estimates that we do not have in our framework.  That result states the uniform convergence of $v(\cdot, t) - \bar\lambda t$ to a critical solution, 
where $v$ is the solution of the Cauchy problem generated by the Lax-Oleinik semigroup. This is not known to be unique in their framework, see Rmk. 9.16 of \cite{CanMen23}. For our quadratic  Lagrangian and $b\equiv 0$ we can combine their result and  our uniqueness Theorem \ref{teo:confronto}  to obtain the uniform convergence of $w(\cdot, t) - \bar\lambda t$ to a critical solution.
\end{rem}
\subsection{Fixed points of the Hamilton-Jacobi semigroup}
\label{sec:fixed}
For any bounded function $\phi\in B(\R^n)$ define the value function of the optimal control problem with finite horizon $t\geq 0$ and terminal cost $\phi$
  \begin{equation}  
\label{sg}
T_t \phi(x):= \inf_{\alpha\in\mathcal{A}}\left\{ \int_0^t  L(y_x(s, \al), \alpha(s))\,ds+\phi(y_x(t, \al))\right\} 
\end{equation} 
where $y_x(\cdot, \al)$ is the trajectory of \eqref{syst}  with $y_x(0, \al)=x$. Note that $T_t \phi\in B(\R^n)$, because $L$ is bounded from above for any fixed control, and it is bounded from below. If $\phi$ is continuous, $w(x, t)=T_t \phi(x)$ is also the unique viscosity solution of the Cauchy problem
\begin{equation} 
\label{Cauchy} 
\partial_t w +H(x, D_xw) = 0 ,\quad \text{in } 
\R^n\times \mathopen( 0,+\infty \mathclose) , \qquad w(\cdot,0)= \phi ,
\end{equation} 
by Thm. \eqref{teo:confronto}. The basic properties of $T_t : B(\R^n) \to B(\R^n)$ that we need are collected in the next proposition.
\begin{prop}
\label{semigroup}
1) For all $\phi\in B(\R^n)$ and $c\in\R$, $T_t (\phi + c) = T_t \phi + c$ ;

\noindent 
2) $\phi\leq\psi$ implies  $T_t \phi \leq T_t \psi$ for all $t\geq 0$ ;


\noindent 
3) $T_t$ has the semigroup property, i.e.,   $T_t(T_s \phi) = T_{t+s} \phi$ for all $\phi \in B(\R^n)$ .
\end{prop}
\begin{proof}
Properties 1) and 2) are immediate consequences of the definition. 
Property 3) expresses the Dynamic Programming Principle for the finite horizon control problem, see, e.g., Prop. III.3.2 in \cite{BCD}. 
\end{proof}
\begin{rem} In the classical weak KAM theory the semigroup is defined with exchanged roles of the initial and the terminal conditions of the trajectories, and it is called 
 Lax-Oleinik semigroup \cite{fathi1997, fathi2008, Fa2014}. As remarked in Sect. 7.1 of \cite{CanMen23}, that definition has a backward character, whereas \eqref{sg} is forward in time and fits better the current context.
\end{rem}
As in  the  classical weak KAM theory we are interested in finding a 
 fixed point $\chi$ of the semigroup $T_t - \bar\lambda t$, which is usually called a {\em weak KAM solution}. This means
  \begin{equation}  
\label{fixed}
\chi(x) + \bar\lambda t = T_t \chi (x) 
 \quad \forall\, x\in \R^n, \,\forall  \,t\geq 0 .
\end{equation}
 Note that any weak KAM solution is also a critical solution, i.e.,  $\chi$ is a non-continuous solution of \eqref{eq:crit5} and $\chi_*$ is a lsc-solution of  \eqref{eq:crit5}, by applying the Dynamic Programming Principle.

The next result says that the converse is also true in two particular cases, i.e., any critical solution, such as those found in Theorems \ref{kam} and \ref{kam2}, is also a weak KAM solution. 
\begin{prop}
\label{fixed point1}
(i)  If the control set $A$ is bounded, then any lsc-solution $\chi\in BLSC(\R^n)$ of $\bar\lambda + H(x, D\chi) = 0$ satisfies also 
\[
\chi(x) + \bar\lambda t = (T_t \chi)_*(x)  \quad \forall x\in \R^n, \forall  t\geq 0 .
\]
\noindent
(ii) If $A =\R^m$ and $\chi\in BC(\R^n)$ is a continuous solution of $\bar\lambda + H(x, D\chi) = 0$, then it satisfies \eqref{fixed}.
\end{prop}
\begin{proof}
Both statements use that the function $\chi(x) + \bar\lambda t$ solves the Cauchy problem \eqref{Cauchy} with initial data $\phi = \chi$. 

In case (ii),  $\chi(x) + \bar\lambda t$ and $T_t \chi$ are continuous in $\R^n \times [0, +\infty)$ and attain continuously the initial data. Then we get the equality from the Comparison Theorem \ref{teo:confronto}.

For case (i), it is known from  \cite{BJ91} that $w(\cdot,t) = (T_t \chi)_*$ is a lsc-solution of \eqref{Cauchy} with $\phi = \chi$ and the additional property that $(w_{|\{t>0\}})_*(x,0)= \chi(x)$, see also \cite{BCD}.
Then the equality \eqref{fixed} follows from the uniqueness theorem for lsc-solution holding under the assumption $|H(x,p)-H(z,p)|\leq C |x-y||p| + \omega(|x-y|)$, where $\omega$ is a modulus, see \cite{BJ, BCD}.
\end{proof}
\begin{rem}
 If $b\equiv 0$ there is a continuous critical solution $\chi$ by Rmk. \ref{contchi}. 
Then {\em (ii)} of the last proposition applies in this case.
\end{rem}
In the general case we can prove 
 the existence of at least one fixed point of the semigroup, and that it is also a critical solution.
\begin{thm}
\label{laxol} 
Assume $H$ is defined by \eqref{Hquad}.
Then there exists $\chi \in B(\R^n)$ satisfying \eqref{fixed}.  Moreover $\chi$ is a non-continuous viscosity solution of the critical equation \eqref{eq:crit2} and 
$\chi_*$ is a lsc-solution of \eqref{eq:crit2}. 
\end{thm}

\begin{proof}
Consider the functions $\underline w \leq \bar w$ found in Thm. \ref{kam2}. Let us first observe that a bilateral supersolution of the critical equation \eqref{eq:crit2}, such as $\underline w$, satisfies a sub-optimality principle, namely,
\begin{equation}  
\label{subsol}
\underline w(x) \leq  \inf_{\alpha\in\mathcal{A}}\left\{ \int_0^t \left( L(y_x(s, \al), \alpha(s)) - \bar\lambda \right)\,ds+ \underline w(y_x(t, \al))\right\} = T_t \underline w(x)  - \bar\lambda t \quad \forall t\geq 0\,,
\end{equation}
by Theorem 3.2 of \cite{sor99}. In the jargon of weak KAM theory this says that $\underline w$ is dominated by $L- \bar\lambda$.

The candidate $\chi$ is the long-time  limit of $T_t \underline w - \bar\lambda t$. We claim 
that this function is non-decreasing in $t$. By the inequality \eqref{subsol} and the properties 
 of the semigroup  in Prop. \ref{semigroup} we have
\[
T_t \underline w (x) - \bar\lambda t \leq T_t (T_s \underline w - \bar\lambda s)(x)  - \bar\lambda t =  T_{t+s} \underline w(x) - \bar\lambda (t+ s) \,,
\]
for all $s>0$, 
which gives  the claim. To find a bound from above we observe that
$\underline w + \bar\lambda t + C$ is a supersolution of the Cauchy problem \eqref{Cauchy} with initial data $\phi =  \bar w$, for $C$ large enough, whereas $T_t \bar w (x)$ is a subsolution.
Then property 2) of Proposition \ref{semigroup} and the Comparison Principle Thm. \ref{teo:confronto} give
\[
T_t \underline w (x)  
\leq  T_t \bar w(x)  \leq \bar C + \bar\lambda t \quad \forall\,x, t> 0 .
\]
Therefore 
\begin{equation}  
\label{chi}
\lim_{t\to +\infty} \left( T_t \underline w(x) - \bar\lambda t \right) = \sup_{t \geq 0} \left( T_t \underline w(x) - \bar\lambda t \right) =: \chi(x)
\end{equation}
defines a function $\chi \in B(\R^n)$. Then the properties of the semigroup in Prop. \ref{semigroup} imply
\begin{multline*}
 T_t \chi(x) - \bar\lambda t =  T_t \left(\sup_{s \geq 0} \left( T_s \underline w - \bar\lambda s \right) \right)(x)- \bar\lambda t =\\ 
\sup_{s \geq 0}  T_t\left( T_s \underline w- \bar\lambda s \right)(x)  - \bar\lambda t =
\sup_{s \geq 0}  \left(T_{t+s} 
 \underline w(x) - \bar\lambda (s + t)\right) = \chi(x) 
\end{multline*}
for all $t>0$, which is the desired fixed point property \eqref{fixed}.

For the last statements we observe that the fixed point equation \eqref{fixed} is a sub- and a superoptimality principle for a control problem with running cost $L-\bar\lambda$. Then it is well known that the upper and lower envelopes of $\chi$ are, respectively, a sub- and a supersolution of \eqref{eq:crit}, see, e.g., Thm. V.2.6 in \cite{BCD}.
Moreover, the suboptimality principle 
\[
\chi(x) \leq  \inf_{\alpha\in\mathcal{A}}\left\{\int_0^t \left( L(y_x(s, \al), \alpha(s)) - \bar\lambda \right)\,ds +\chi(y_x(t, \al))\right\} , \quad \forall t\geq 0\,,
\]
implies the backward Dynamic Programming principle
\begin{equation*}  
\chi(x) \geq  \chi(y_x(-t, \al)) - \int_0^t \left( L(y_x(-s, \al), \alpha(-s)) - \bar\lambda \right)\,ds , \quad \forall t\geq 0\,,
\end{equation*}
for all  $\al : (-\infty, 0] \to A$ measurable and $x\in\R^n$, as in 
Lemma V.5.4 of  \cite{BCD}. Then 
 $\chi_*$ is also a supersolution of $-\bar\lambda-H(x, D\chi_*)=0$, as in Prop. V.5.2 of  \cite{BCD}.
\end{proof}

\section{Application to periodic homogenization}
\label{homo}
In this section we consider the problem of the convergence as $\ep\to0+$ of the solutions of the equation 
\begin{equation}  
\label{pdehomo}
v ^\ep+ \frac 12 \left|\sigma\left(\frac{z}{\ep}\right)D_zv^\varepsilon\right|^2+f\left(\frac{z}{\epsilon}\right)\cdot D_zv^\epsilon = g\left(z,\frac{z}{\epsilon}\right) ,\quad z\in\R^n
\end{equation}
and of its evolutive counterpart, where $\sigma$ is as in Section \ref{Prelim} with the vector field $F$ satisfying the bracket condition \eqref{SBG}, $f: \R^n \to\R^n$ is $\mathbb{Z}^n$-periodic and Lipschitz,
and $g \in BUC(\R^n\times\R^n)$ is periodic in the second variable. 
The existence, for each $\ep>0$, of a bounded viscosity solution $v^\ep$ follows form standard method (e.g., Perron-Ishii, or representation as value function of an infinite horizon discounted control problem), and uniqueness comes from Thm. \ref{teo:comp}.

 We want to find an effective  Hamiltonian $\bar H$ such that the limit of $v^\ep(z)$ satisfies
\begin{equation*} 
v +\bar H(z,D_zv) = 0 .
\end{equation*}
To guess who $\bar H$ can be we make the formal expansion $v^\ep(z)= v(z) +\ep u(\frac{z}{\epsilon}) + \text{ l.o.t.}$ and plug it into the equation to get
\[
v +  \frac 12 |\sigma\left(\frac{z}{\ep}\right) (D
v + Du\left(\frac{z}{\ep}\right))|^2+f\left(\frac{z}{\ep}\right)\cdot (D
v + Du\left(\frac{z}{\ep}\right)) = g\left(z,\frac{z}{\ep}\right) + \text{ l.o.t.} .
\]
This suggests to freeze $z$ and $P=Dv(z)$, set $x=z/\ep$, and consider 
$
\bar H(z,P) = -\bar \lambda ,
$
where $\bar \lambda$ and $u(x)$ solve the {\em cell problem}
\[
\bar\lambda +  \frac12 |\sigma(x)(P+Du)|^2 + f(x)\cdot(P+ Du) = g(z, x) \quad 
x\in \R^n .
\]
This becomes the critical equation \eqref{eq:crit2}, with  parameters $z$ and $P$, once we set
\begin{equation} 
\label{data}
b(x):= f(x) + \sigma(x)^T\sigma(x) P ,\quad l(x) := g(z, x) -  f(x)\cdot P -  \frac 12 |\sigma(x) P|^2 .
\end{equation}
We define the {\em effective Hamiltonian} by the critical value defined in Thm. \ref{kam} for the data $b, l$ given by \eqref{data}
\begin{equation} 
\label{effH}
\bar H(z,P) := - \bar\lambda \,.
\end{equation}
\begin{thm}
\label{convhomo}
Under the assumptions listed above, the solution $v^\ep(z)$ of \eqref{pdehomo} converges locally uniformly 
 as $\ep\to 0+$ to the unique bounded solution of 
 \begin{equation} 
\label{eff}
v +\bar H(z,D
v) = 0 \;\text{ in } \R^n .\end{equation}
\end{thm}
\begin{proof}
We follow the approach to homogenization via singular perturbations of \cite{ABsiam, ABarma, ABmem}.  We double the space variables and consider the PDE
\begin{equation} 
\label{sp}
u^\epsilon + \frac 12 \left|\sigma(x)\left(D_zu^\epsilon + \frac{D_xu^\epsilon}\ep\right) \right|^2+f(x)\cdot \left(D_zu^\epsilon + \frac{D_xu^\epsilon}\ep\right) = g(z,x) ,\quad  
\text{in }\R^{n}\times \R^{n} .\end{equation}
This equation  has a  bounded viscosity solution $u^\ep(z,x)$ by standard methods, unique by Thm. \ref{teo:comp}, and it is easy to check that
\[
v^\ep(z)=u^\ep(z,\frac z\ep) \,,
\]
where $v^\ep$ solves \eqref{pdehomo}. The Hamiltonian in \eqref{sp} fits the general conditions of \cite{ABarma, ABmem} except the regularity in $x$ required  there for using standard comparison principles. In  \eqref{sp} the Hamiltonian is quadratic, but we have the comparison principles stated in Sect. \ref{Prelim}. 
Moreover such Hamiltonian is ergodic in the sense of \cite{ABarma}. To see it 
consider the associated discounted equation with $z$ and $P$ frozen
\begin{equation} 
\label{cellpz}
\delta w_\de +  \frac12 |\sigma(x)(P+Dw_\de )|^2 + f(x)\cdot(P+ Dw_\de ) = g(z, x)  \qquad 
x\in\R^n ,
\end{equation}
and use \eqref{data} to rewrite it as
\[
\delta w_\de +  \frac12 |\sigma(x)(Dw_\de )|^2 + b(x)\cdot Dw_\de  = l(x)  \qquad 
x\in\R^n ,
\]
where $b$ and $l$ depend on the parameters $z, P$. Then Corollary \ref{conv2} gives
\begin{equation} 
\label{limit}
\lim_{\de\to 0} \de w_\de(x) = -\bar H(z,P) \qquad\text{uniformly in } \R^n ,
\end{equation}
This says that $w_\de$ can be used as an approximate corrector in the perturbed test function method adapted from Evans \cite{E92} to singular perturbations \cite{ABarma}. This theory entails that, if $u^\ep$ are locally equibounded, the semilimits 
\[
\bar u(z) := \limsup_{\ep\to 0, z'\to z} \sup_x u^\ep(z',x) ,\quad \underline u(z) := \liminf_{\ep\to 0, z'\to z} \inf_x u^\ep(z',x)
\]
are,  respectively, a sub- and a supersolution of the effective PDE  \eqref{eff}, see Thm. 2 in \cite{ABarma}. Since $\pm \sup |g|$ are a super- and a subsolution of \eqref{sp}, the comparison principle gives the estimate 
\[
-\sup |g|  \leq u^\ep(z, x)\leq \sup |g| 
\]
and therefore the equiboundedness.

We want to show that the effective equation \eqref{eff} satisfies  the comparison principle. The continuity of $(z,p)\mapsto \bar H(z,p)$ is a known consequence of ergodicity, see e.g. Prop. 3 in \cite{ABarma}. We claim that 
\begin{equation} 
\label{estH}
|\bar H(z,p) - \bar H(z',p)|\leq \omega(|z-z'|) \quad \forall \, z, z', p\in \R^n ,
\end{equation}
where $\omega$ is the modulus of continuity of $g$, and therefore of $l$, with respect to $z$. To prove this estimate we compare the solution $w_\de$ of \eqref{cellpz} with the solution $w'_\de$ of the same equation with $z$ replaced by $z'$. Then $w'_\de$ satisfies 
\[
\delta w'_\de +  \frac12 |\sigma(x)(Dw'_\de )|^2 + b(x)\cdot Dw'_\de  \leq l(z, x)  + \omega(|z-z'|)
\]
and so $w'_\de-\frac 1\de  \omega(|z-z'|)$ is a subsolution of of \eqref{cellpz}. Then the comparison principle Thm. \ref{teo:comp} gives
\[
w'_\de \leq w_\de + \frac 1\de  \omega(|z-z'|) 
\]
and by \eqref{limit} we get 
\[
\bar H(z,p) - \bar H(z',p) \leq \omega(|z-z'|) . 
\]
By inverting the roles of $z$ and $z'$ we get \eqref{estH}. Then the comparison principle for \eqref{eff} is a standard result (see, e.g., \cite{BCD}) . We use it for the semilimits and get that $\bar u \leq \underline u$, and so $\bar u = \underline u =: v$ is the unique bounded solution of \eqref{eff}.

Finally we go back to the solution $v^\ep$ of \eqref{pdehomo}. Clearly 
\[
\bar v(z):= \limsup_{\ep\to 0, z'\to z}  v^\ep(z')  \leq \bar u(z) , \qquad \underline v(z) := \liminf_{\ep\to 0, z'\to z} v^\ep(z') \geq\underline u(z) .
\]
Then $\bar v = \underline v = v$, which implies the convergence of $v^\ep$ to $v$ locally uniformly.
\end{proof}
\begin{rem}
\label{homAbdd}
In the same way we can prove the locally uniform convergence of the solutions of the equation 
\begin{equation*}  
v ^\ep+   b\left(\frac{z}{\epsilon}\right)\cdot Dv^\epsilon +\sup_{a\in A} \left\{ Dv^\ep \cdot F\left(\frac{z}{\ep}\right) a - l_1\left(\frac{z}{\ep}, a\right) \right\} = g\left(z,\frac{z}{\epsilon}\right) ,\quad z\in\R^n 
\end{equation*}
where $A$ is bounded and contains $B_m(K)$. The cell problem in this case is
\[
\bar\lambda +   b(x)\cdot(P+ Du)+ \sup_{a\in A} \left\{ (P+ Du) \cdot F(x) a - l_1(x, a) \right\}  = g(z, x) 
\quad 
x\in \R^n ,
\]
which can be written in the form of the critical equation \eqref{eq:crit5}
\[
\bar\lambda +   b(x)\cdot Du+ \sup_{a\in A} \left\{ Du \cdot F(x) a - L(x, a, z, P) \right\} ,
\]
where $z$ and $P$ are parameters and
\[
L(x, a, z, P):= l_1(x, a) + g(z, x) - P\cdot  (b(x) +F(x) a).
\]
We set again $\bar H(z,P)=-\bar \la$ and obtain that the limit $v$ of $v ^\ep$ is the unique solution of the effective equation \eqref{eff}.
\end{rem}

Finally we consider  the homogenization of solutions of the Cauchy problem with oscillating initial data
\begin{equation*}  
\begin{sistema}
v_t^\epsilon + \frac 12 |\sigma(\frac{z}{\epsilon})D_zv^\epsilon|^2+f(\frac{z}{\epsilon})\cdot D_zv^\epsilon = g(z,\frac{z}{\epsilon}) ,\quad (z,t)\in\R^n\times \mathopen( 0,T \mathclose) ,\\
\\
v^\epsilon(z,0)=h(z,\frac{z}{\epsilon}),
\end{sistema}
\end{equation*}
where  $h\in BUC(\R^n)$ is periodic.
The existence, for each $\ep>0$, of a bounded viscosity solution is again standard, and uniqueness comes from Thm. \ref{teo:confronto}.

\begin{thm}
\label{convhomo2}
Under the assumptions listed above, the solution $v^\ep(z,t)$ converges locally uniformly in $\R^n\times (0, T)$ as $\ep\to 0$ to the unique bounded solution of the Cauchy problem
\begin{equation} 
\label{effCauchy}
v_t +\bar H(z,D_zv) = 0 \;\text{ in } \R^n\times \mathopen( 0,T \mathclose) ,\quad  v(z,0) = \bar h(z) := \min_x h(z, x) , \;\forall z ,
\end{equation}
where $\bar H$ is defined by \eqref{effH}
\end{thm}
\begin{proof} Most of the proof is the same as that of Thm. \ref{convhomo}, so we only discuss the differences. We consider the evolutive version of the singularly perturbed equation \eqref{sp} where $u^\ep$ is replaced by $u^\ep_t$ and with the initial condition $u^\ep(z,x,0) = h(z, x)$. Then
\[
v^\ep(z,t)=u^\ep(z,\frac z\ep,t) \,, \quad \text{in } \R^n\times [0,T \mathclose) .
\]
We have  the estimate 
\[
-t\sup |g| - \sup |h| \leq u^\ep(z, x, t)\leq t\sup |g| + \sup |h| ,
\]
and the semilimits, for $t>0$,  
\[
\bar u(z,t) := \limsup_{\ep\to 0, (z',t')\to (z,t)} \sup_x u^\ep(z',x,t') ,\quad \underline u(z,t) := \liminf_{\ep\to 0, (z',t')\to (z,t)} \inf_x u^\ep(z',x,t')
\]
are 
 a sub- and a supersolution of the effective PDE in \eqref{effCauchy}, by Thm. 2 in \cite{ABarma}. 
 
 It remains to handle the boundary layer at the initial time. We extend $\bar u$  and $\underline u$ by
 \[
\bar u(z,0) := \limsup_{(z',t')\to (z,0+)} \bar u(z', t') ,\qquad \underline u(z,0) := \liminf_{(z',t')\to (z,0+)} \bar u(z', t') .
\]
The initial condition satisfied by these semilimits comes from the property of stabilization to a constant of the solution to the homogeneous equation 
\begin{equation} 
\label{stabil}
w_t + \frac12 |\sigma(x)D_xw|^2 = 0 \;\text{ in } \R^n\times \mathopen( 0,+\infty\mathclose) ,\quad  w(x,0) =  h(z, x) , \;\forall z ,
\end{equation}
where $z$ is a frozen parameter, namely
\[
\lim_{t\to +\infty} w(x, t; z) = \text{const.} =: \bar h(z) .
\]
Then Thm. 3 in \cite{ABarma} gives 
\[
\bar u(z,0) \leq \bar h(z) \quad\text{and}\quad  \underline u(z,0) \geq \bar h(z) ,
\]
and so $\bar u(\cdot,0) =\underline u(\cdot,0) =  \bar h$. (Actually, in \cite{ABarma} the Hamiltonian is 1-homogeneous and here it  is 2-homogeneous, but the proof holds without changes).
The stabilization property of the equation \eqref{stabil} comes from the fact that it is the HJB equation of an optimal control problem for  the driftless control system 
$\dot{y}(t)=  - F(y(t)) \al(t)$. If $A\supseteq B_m(K)$ 
 this system satisfies \eqref{BTC}, by the assumption \eqref{SBG}, and moreover it is stoppable, i.e.,  for each point of the state space  there is a control allowing to remain there forever. Then Prop. 6.5 in \cite{ABmem} gives the conclusion, as well as the explicit formula $\bar h(z) = \min_x h(z, x)$.

 Now the comparison principle for the Cauchy problem \eqref{effCauchy}, which holds by the estimate \eqref{estH} on $\bar H$, implies that $\bar u = \underline u =: v$ is the unique bounded solution of \eqref{effCauchy}. Finally, as in the previous proof, we consider the weak limits of $v^\ep$, compare them with $\bar u$ and $\underline u$, and conclude the uniform convergence of  $v^\ep$ locally uniformly in $\R^n\times(0,T)$.
\end{proof}
\begin{rem}

As in Remark \ref{homAbdd} we can consider also the case of bounded $A\supseteq B_m(K)$ and the 
Cauchy problem 
\begin{equation*}  
\begin{sistema}
v_t^\epsilon +  b\left(\frac{z}{\epsilon}\right)\cdot Dv^\epsilon +\sup_{a\in A} \left\{ Dv^\ep \cdot F\left(\frac{z}{\ep}\right) a - l_1\left(\frac{z}{\ep}, a\right) \right\} = g\left(z,\frac{z}{\epsilon}\right) ,\quad (z,t)\in\R^n\times \mathopen( 0,T \mathclose) ,\\
\\
v^\epsilon(z,0)=h(z,\frac{z}{\epsilon}).
\end{sistema}
\end{equation*}
By the same proof we obtain again that the relaxed semilimits $\bar u$ and $\underline u$ are a sub- and a supersolution of the effective equation $v_t +\bar H(z,D_zv) = 0$ with $\bar H$ defined in Remark \ref{homAbdd}. We obtain also the locally uniform convergence of $v^\ep$ to the unique solution of \eqref{effCauchy} if either $b\equiv 0$ or the initial data do not oscillate, i.e., $h=h(z)$. In the general case the existence of some effective initial data $\bar h$ does not follows from known results. In fact, the homogeneous equation arising in the study of stabilization is not \eqref{stabil}, it is instead 1- homogeneous in $D_xw$ and has the term $b\cdot D_xw$. Hence the associated control system is not driftless and therefore not stoppable in general. This case is left open to future research.
\end{rem}

\end{document}